\renewcommand{\to}{\rightarrow}

\documentclass[11pt]{article}

\usepackage{amsfonts, amssymb, amsmath, graphicx}
\usepackage[all, knot]{xy}

\newcommand{\N}{\mathbb{N}}

\newcommand{\Z}{\mathbb{Z}}
\newcommand{\R}{\mathbb{R}}

\def\calA{{\mathcal A}}
\def\calC{{\mathcal C}}
\def\calH{{\mathcal H}}
\def\calF{{\mathcal F}}
\def\Kof{{\textit Kof}}

\newtheorem{theorem}{Theorem}
\newtheorem{corollary}[theorem]{Corollary}
\newtheorem{proposition}[theorem]{Proposition}

\input{epsf.sty}
\input xy
\xyoption{all}

\def\co{\colon\thinspace} 
\DeclareMathOperator{\rk}{rk}

\newcommand{\cat}[1]{\textit{#1}}
\newcommand{\rmod}{\cat{R-mod}}

\newcommand{\cob}{\cat{Cob}}
\newcommand{\foams}{\cat{Foams}}

\begin{document}
\title{A cohomology theory for colored tangles}

\author{
Carmen Caprau\footnote{The author was supported in part by NSF grant DMS 0906401}
\\ California State University, Fresno}

\maketitle

\begin{abstract}
We employ the $sl(2)$ foam cohomology to define a cohomology theory for oriented framed tangles whose components are labelled by irreducible representations of $U_q(sl(2))$. We show that the corresponding colored invariants of tangles can be assembled into invariants of bigger tangles.  For the case of knots and links, the corresponding theory is a categorification of the colored Jones polynomial, and provides a tool for efficient computations of the resulting colored invariant of knots and links. Our theory is defined over the Gaussian integers $\Z[i]$ (and more generally over $\Z[i][a,h]$, where $a,h$ are formal parameters), and enhances the existing categorifications of the colored Jones polynomial.
\vspace{1cm}


\end{abstract}

\section{Introduction} \label{introduction}

In~\cite{Kh3} Khovanov constructed a cochain complex associated to an oriented framed link whose components are labelled by irreducible representations of $U_q(sl(2))$. The graded Euler characteristic of the homology of this complex is the colored Jones polynomial. Specifically,~\cite{Kh3} provides a categorification of the colored Jones polynomial by interpreting the defining formula for the polynomial 

\[ J_n(K) = \sum_{i = 0}^{\lfloor \frac{n}{2} \rfloor}(-1)^i\binom {n-i}{i} J(K^{n-2i}),\]
where $K^j$ is the $j$-parallel cable of the knot $K$, as the Euler characteristic of a complex whose objects require the original Khovanov homology~\cite{Kh1} of the cablings $K^{n-2i}$, for $i = 0, \dots, \lfloor \frac{n}{2} \rfloor$. As a consequence of the fact that the original Khovanov homology is functorial with respect to link cobordisms only up to a negative sign, the construction in~\cite{Kh3} works over $\Z_2$.  

Mackaay and Turner~\cite{MT} followed Khovanov's proposed categorification~\cite{Kh3} of the colored Jones polynomial to define and compute the colored Bar-Natan theory. Mackaay-Turner work uses Bar-Natan's~\cite{BN} filtered theory in place of Khovanov's theory,  and is defined, once again, over $\Z_2$.  

In~\cite{BW} Beliakova and Wehrli developed homology theories of colored links over $\Z[1/2]$ by using Bar-Natan's geometric \textit{formal Khovanov bracket}, which is an object in the category Kob:= Kom(Mat($\cob^3_{/\ell}$)). Here $\cob^3_{/\ell}$ is the category of $2$-cobordisms modulo some local relations $\ell$. For any additive category $\calC$, Mat($\calC$) is the category whose objects are formal direct sums of objects of $\calC$, and whose morphisms are matrices of morphisms in  $\calC$ for which the composition law is modeled on matrix multiplication. Kom($\calC$) is the category of chain complexes over $\calC$ whose objects are bounded (co)chain complexes, and whose morphisms are chain transformations. The colored link invariant in~\cite{BW} is a complex whose objects are formal direct sums of formal Khovanov brackets. Beliakova and Wehrli explained that there is a way to deal with the sign ambiguity in the functoriality property of the formal Khovanov bracket, without the need of working over a field of characteristic two (they showed that there is a \textit{satisfactory choice of signs} making all squares in the cube of resolutions associated to a link diagram  anticommutative). Their arguments imply that the categorification defined in~\cite{Kh3} works over integers. 

The goal of this paper is two-folded. First we enhance the existing categorifications of the colored Jones polynomial by giving a clean definition of the colored invariant of a knot or a link, in the sense that one is not restricted to work over $\Z_2$. For that, we employ the universal $sl(2)$ cohomology theory that uses foams (also called seamed cobordisms, or singular cobordisms), constructed by the author in~\cite{C2} (compare with~\cite{C1}), and which is properly functorial under link cobordisms with no sign ambiguity. The theory developed in~\cite{C2} is a Khovanov-type tangle cohomology theory  defined over the ring $R = \Z[i][a,h]$, where $i^2 = -1$, and $a$ and $h$ are formal parameters. We refer to this theory as the (universal) \textit{$sl(2)$ foam cohomology}. We won't make any specific computations here, so one can let either $a$ or $h$ be zero (or both, for that matter). Cabling a knot or a link introduces an unmanageable number of crossings from a computational point of view, which brings us to the second goal of the paper, namely to define a theory in which the invariants can be computed efficiently. For that, we use  Bar-Natan's ``divide and conquer" approach to computations and construct a \textit{local} colored cohomology theory, in that it is built with colored tangles in mind and which composes well under tangle composition. 

We construct a triply-graded cohomology theory for colored oriented framed tangles, which for the case of links (that can be considered as closed tangles) is a categorification of the colored Jones polynomial. The resulting invariant of tangles has excellent composition properties, allowing one to obtain the invariant of a colored framed link from the invariants of its subtangles. We will discuss the functoriality property of our invariant with respect to tangle cobordisms (rel. boundary) in a subsequent paper.

The paper is organized as follows. Section~\ref{review} overviews the main facts about the universal $sl(2)$ foam cohomology which will be extensively used in this paper, and recalls the definition of the colored Jones polynomial of an oriented framed link. In Section~\ref{sec:colored link-coho} we define the new cohomology theory for colored framed links, categorifying the colored Jones polynomial. Section~\ref{sec:colored tangle-coho} is dedicated to tangles. We first consider the case of mono-colored tangles with no closed components, and construct a cohomology theory for such tangles. Then, we generalize our construction to arbitrary colored framed tangles. In both cases, we show that there is a composition rule that takes the invariants of tangles to invariants of bigger tangles, and thus produces the invariant of a knot or a link.

\section{Brief review of necessary concepts}\label{review}
\subsection{Universal $sl(2)$ foam cohomology}

We assume familiarity with the construction in~\cite{C2}, but we briefly recall a few concepts, notations and results coming from the (universal) $sl(2)$ foam cohomology theory, which are essential in understanding this paper. The construction involves \textit{webs} and dotted \textit{foams} modulo local relations, along the lines of~\cite{BN} and~\cite{Kh2}. A web here is a planar graph with bivalent vertices which are either ``sources" or ``sinks", and a  foam is a 2-cobordism between webs. Foams contain \textit{seams}, where a seam is a singular arc or a singular circle where orientations disagree. The author constructed in~\cite{C2} a doubly graded cohomology theory (for oriented tangles) over the graded ring $R = \Z[i][a,h],$ where $i^2 = -1$, and $a$ and $h$ are parameters with $\deg(a) = 4, \deg(h) = 2$. (The special case of $h = 0$ was treated in~\cite{C1}; see also its longer and more detailed preprint version~\cite{C0}.)

We denote by $\foams$ the additive category whose objects are webs and whose morphisms are $R$-linear combinations of dotted foams, and we denote by $\foams_{/\ell}$ the quotient category of $\foams$ by a finite set of relations $\ell$; that is, we mod out the morphisms of the category $\foams$ by the local relations $\ell$ -- these are generalized Bar-Natan relations~\cite{BN} for the local Khovanov homology, enhanced by additional relations involving the 2-sphere with a seam. There is a functor $\calF \co \foams_{/\ell} \to \rmod$ taking us from the geometric picture to the algebraic picture, where $\rmod$ is the category of $R$-modules and module homomorphisms. This functor is related to the \textit{universal  rank-two Frobenius system} defined on the $R$-module $\calA = R[X]/(X^2 - hX -a)$, graded by $\deg(1) = -1, \deg(X) = 1$ (the universal Frobenius system was coined by Khovanov in~\cite{Kh4}). With respect to the generators $1$ and $X$ of the algebra $\calA$, the counit and comultiplication maps  are given by $\epsilon(1) = 0, \epsilon(X) = 1$ and $\Delta(1)= 1\otimes X + X \otimes 1 - h 1 \otimes 1, \Delta(X) = X \otimes X + a 1 \otimes 1$, respectively. A dot on a foam corresponds to the multiplication by $X$ endomorphism of $\calA$.
 The TQFT corresponding to $\mathcal{A}$ factors through the quotient category of \textit{Foams} by the relations $\ell.$
  
  We denote by $\Kof = $ Kom(Mat($\foams_{/\ell}))$ the category of complexes whose objects are column vectors of webs and whose morphisms are matrices of dotted foams modulo the local relations $\ell$. Moreover, let $\Kof_{/h} = $ Kom$_{/h}$(Mat($\foams_{/\ell}))$ be the homotopy subcategory of the earlier. 
    
  Given a planar diagram $D$ of an oriented tangle $T$, we constructed in~\cite{C2} a formal cochain complex $[D] \in \Kof$, whose homotopy class is an invariant of $T$. By applying a degree-preserving Bar-Natan type functor $\calF \co \foams_{/\ell} \to \rmod$, which extends to a functor $\calF \co \Kof \to \mbox{Kom(Mat}(\rmod))$, we obtain a cochain complex $\calF[D]$ whose homology $\calH(D): = H(\calF[D])$ is a doubly-graded invariant of $T$. If the tangle $T$ is a link, the graded Euler characteristic of the homology group $\calH(D)$ is the quantum $sl(2)$ polynomial of the link. Finally, we define $[T]: = [D]$, for any planar diagram $D$ representing the tangle $T$.
   
  \begin{proposition}\label{prop:functoriality}
Let $C \subset \R^3 \times [0,1]$ be a tangle cobordism between tangles $T_1$ and $ T_2$, and denote by $B$ the set of boundary points of $T_1$ (and thus of $T_2$). There exists an induced graded map $[T_1] \to [T_2]$ of degree $-\chi(C) + \frac{1}{2}|B|$, well-defined under ambient isotopy of $C$ (rel. boundary), where $\chi(C)$ is the Euler characteristic of $C$, and $|B|$ is the cardinality of $B$.
\end{proposition}

For the scope of this paper, it is important to recall that the geometric invariant has excellent composition properties, making the $sl(2)$ foam cohomology theory ready for a ``divide and conquer" approach to computations. One can cut a link $L$ into subtangles $T_i$, compute the geometric invariant $[T_i]$ for each of these tangles, and finally assemble the obtained invariants into the invariant of $L$ via the tensor product operation induced on formal complexes. But before performing the assembling operation one can simplify each $[T_i]$ as much as possible, via the \textit{delooping} and \textit{Gaussian elimination} procedures. These techniques provide computational efficiency of the $sl(2)$ foam cohomology groups (and, implicitly, of the original Khovanov homology groups). For more details about efficient computations we refer the reader to~\cite{C3}. In particular, it follows that the category $\Kof_{/h}$ has a natural structure of an oriented planar algebra.

\begin{proposition} \label{prop:planar-algebra}
$[ \,\cdot \, ]$ is a planar algebra morphism from the planar algebra of oriented tangles modulo the three Reidemeister moves to the planar algebra $\Kof_{/h}$. 
\end{proposition}

The categories $\foams$ and $\foams_{/\ell}$ are \textit{canopolies} over the planar algebra of web diagrams.  A canopoly is both a category and a planar algebra. This term was coined by Bar-Natan in~\cite{BN}. The category $\Kof$ (and hence $\Kof_{/h}$) can also be viewed as a canopoly, where the ``tops" and ``bottoms" of  cans are formal complexes over $\foams_{/\ell}$, and the ``cans" are morphisms between complexes. Cobordisms between oriented tangle diagrams can also be composed like tangles, by placing them next to each other and connecting the common ends. Therefore, they form a planar algebra, and thus the category $\cob^4$ of cobordisms between oriented tangle diagrams is a canopoly over the planar algebra of oriented tangle diagrams.  

\begin{proposition} \label{prop:canopoly}
$[ \,\cdot \, ]$ descends to a degree preserving canopoly morphism $[ \,\cdot \, ] \co \cob^4_{/i} \to \Kof_{/h}$ from the canopoly of movie presentations of cobordisms between oriented tangle diagrams, up to isotopy, to the canopoly $\Kof_{/h}$ of formal complexes and morphisms between them, up to homotopy.
\end{proposition} 

The proofs of the above three propositions can be found in~\cite{C0}.

\subsection{Colored Jones polynomial} 

Let $\underline{n} = (n_1,  \ldots, n_l)$ be a vector whose entries are non-negative integers. Let $(L, {\underline{n}})$ be an oriented framed link with $l$ components colored by $\underline{n}$; that is, the $i$-th component of $L$ is colored (or labeled) by  $n_i$, or equivalently, by the $(n_i +1)$-dimensional irreducible representation of the quantum group $U_q(sl(2))$. We denote the colored Jones polynomial of $(L, {\underline{n}})$ by $J_{\underline{n}}(L)$. It is a Laurent polynomial in $q$, and is given by the formula

\[ J_{\underline{n}}(L) = \sum_{\underline{k} = \underline{0}}^{\left \lfloor  \underline{n}/2 \right \rfloor} (-1)^{|\underline{k}|} \binom{\underline{n} - \underline{k}} {\underline{k}}  J(L^{\underline{n}-2\underline{k}}),\]
where
\[\ |\underline{k}| = \sum_i k_i \  \ \text{and} \  \   \binom{\underline{n} - \underline{k}} {\underline{k}} = \displaystyle \prod_{i = 1} ^l \binom{n_i - k_i} {k_i}.\] 
Here, $J(L^{\underline{n} - 2 \underline{k}})$ stands for the original Jones polynomial of the $(\underline{n}-2 \underline{k})$-parallel cable of $L$, formed by taking the $(n_i - 2k_i)$-parallel cable of the $i$-th component of $L$ with respect to its framing, for all $1 \leq i \leq l$. 

If all components of $L$ are labeled by $1$, the invariant is the original Jones polynomial of $L$. 

When forming the $m$-parallel cable of a component $K_i$ of $L$, we enumerate the strands in a cross-section of the planar projection of the cable $K_i^m$ from left to right by $1$ to $m$, and orient the parallel cable-strands such that adjacent strands receive opposite orientations, where we give strand $1$ the original orientation of $K_i$.

\section{Colored link cohomology}\label{sec:colored link-coho}

In this section we borrow definitions and some clever ideas from~\cite{BW} and~\cite{Kh3} to construct a cohomology theory for colored framed links, but instead of using the Khovanov homology we employ the universal $sl(2)$ foam cohomology. Although our construction for the case of links is similar to that in~\cite{BW}, it is much simpler.  

Let $(L, \underline{n})$ be an oriented framed link colored by $\underline{n} = (n_1,  \ldots, n_l)\in \N^l$, and let $D$ be a planar diagram for $L$ whose blackboard framing agrees with the given framing of $L$. Denote by $D_1, \dots, D_l$ the components of $D$.

The binomial coefficient $\displaystyle \binom{n - k}{k}$ equals the number of ways of selecting $k$ pairs of neighbors from $n$ dots placed on a line, such that each dot appears in at most one pair. A \textit{dot-row} $s$ is a set of $n$ dots on a line in which some adjacent dots are paired. Denote by $p(s)$ the number of pairs in $s$. Similarly, $ \displaystyle \binom {\underline{n} - \underline{k}} {\underline {k}}$ is the number of ways of selecting $\underline{k}$ pairs of neighbors from $\underline{n}$ dots placed on $l$ lines, where the $i$-th line contains $n_i$ dots. Denote by $\underline{s} = (s_1, \ldots, s_l)$ a set of dot-rows $s_i$ with $n_i$ dots, respectively, and call it a \textit{dot-row vector}. Let $\underline{p}(\underline{s}) = (p(s_1), \ldots, p(s_l))$, and denote by $|\underline{p}(\underline{s})| = p(s_1) + \ldots + p(s_l)$.

\begin{figure}[ht]
\begin{center}
\includegraphics[height=1.7in, angle=90]{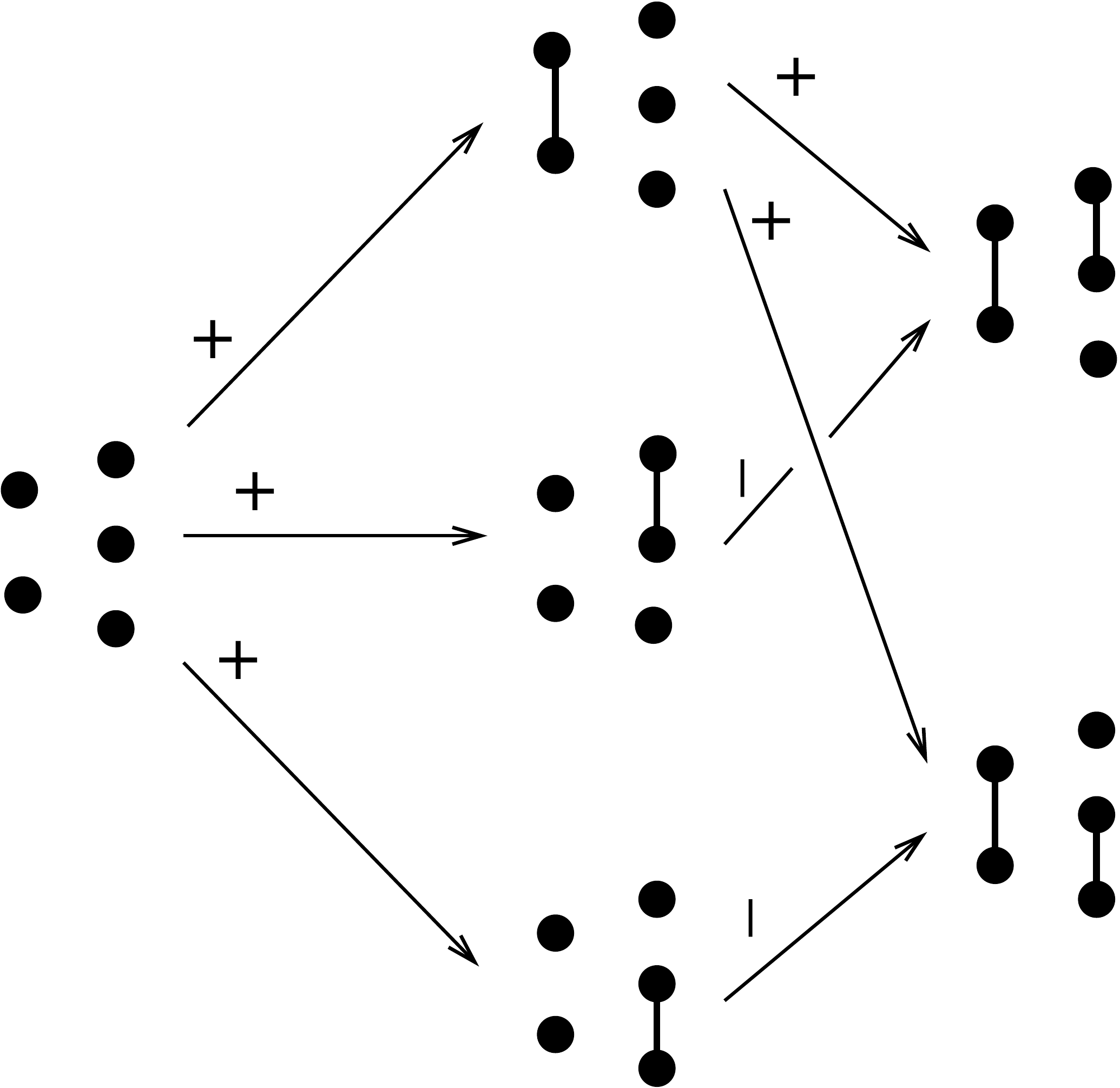} \\
\end{center}
\caption{The graph $\Gamma_{(2,3)}$}\label{fig:graph for a link}
\end{figure}

Let $\Gamma_{\underline{n}}$ be the directed graph whose vertices are in bijective correspondence with all possible dot-row vectors $\underline{s}$ corresponding to $\underline{n}$. Two vertices $\underline{s}$ and $\underline{s}'$ in $\Gamma_{\underline{n}}$ are connected by an edge $e \co \underline{s} \to \underline{s}'$ if and only if all pairs in $\underline{s}$ are pairs in $\underline{s}'$, and $|\underline{p}(\underline{s}')| = |\underline{p}(\underline{s})| + 1$. The \textit{height of a vertex} $\underline{s}$ is equal to $\underline{p}(\underline{s})$, and the edges are oriented from lower to higher heights. In Figure~\ref{fig:graph for a link} we show such a graph for a link with two components colored by $\underline{n}= (2,3)$. 

To a dot-row vector $\underline{s}$ we associate the cable diagram $D_{\underline{s}}: = D^{\underline{n} - 2\underline{p}(\underline{s})}$ formed by taking the $(n_i - 2p(s_i))$-parallel cable of the $i$-th component $D_i$ of $D$, for $1 \leq i \leq l$. In other words, there is a cable-strand for each unpaired dot in $\underline{s}$. To an edge $e \co \underline{s} \to \underline{s}'$ we associate the cobordism $S_e: D_{\underline{s}} \to D_{\underline{s}'}$ given by contracting the neighboring strands in $D_{\underline{s}}$ corresponding to the pair of dots in $\underline{s}'$ but not in $\underline{s}$. In addition, the cobordism $S_e$ receives the sign $(-1)^{o(\underline{s}, \, \underline{s}')}$, where $o(\underline{s},\underline{s}')$ is the number of pairs in $\underline{s}$ to the right and above of the only pair in $\underline{s}' \setminus \underline{s}$. Figure~\ref{fig:graph for a link} explains the sign $(-1)^{o(\underline{s}, \, \underline{s}')}$ and Figure~\ref{fig:edge-map} displays a cobordism $S_e$. The Euler characteristic of the cobordism $S_e$ is $0$, and thus $\mbox{deg}(S_e) = 0$ (we are using the degree convention of the $sl(2)$ foam cohomology theory). 

The cochain complex $C_{\underline{n}}(D)$ for the colored link theory is obtained by applying to the graph $\Gamma_{\underline{n}}$ (with link diagrams as vertices and link cobordisms as oriented edges) the morphism $[ \, \cdot  \,]$ constructed in the $sl(2)$ foam cohomology theory. 
Specifically,
\begin{eqnarray*} D_{\underline{s}} \,\, &\mapsto& \,\,  [D_{\underline{s}}] \in \mbox{Obj}(\Kof_{/h}) \\
D_{\underline{s}} \stackrel{S_e}{\longrightarrow} D_{\underline{s}'}  \,\, &\mapsto& \,\, [D_{\underline{s}}]\stackrel{[S_e]}{\longrightarrow}[D_{\underline{s}'}] \in \mbox{Mor}(\Kof_{/h}).\end{eqnarray*}

According to Proposition~\ref{prop:functoriality}, $[S_e]$ is a well-defined homotopy class of chain maps. The $i$-th cochain object of $C_{\underline{n}}(D)$ is a formal direct sum of complexes at height $i$:
\[ C^i_{\underline{n}}(D): =  \bigoplus_{\underline{s}} [D_{\underline{s}}].\]
The sum above is taken over all dot-row vectors $\underline{s}$ (vertices in $\Gamma_{\underline{n}}$) such that $|\underline{p}(\underline{s})| = i$. The $i$-th differential $d^i : C^i_{\underline{n}}(D) \to C^{i+1}_{\underline{n}}(D)$ is a formal sum of all morphisms $[S_e]$ corresponding to edges $e$ at height $i$. That is, if $v \in [D_{\underline{s}}]$ then $d^i(v): = \sum_{e} [S_e](v)$, where the sum is over all edges $e$ with tail $\underline{s}$. 

Observe that the maps $d^i$ are degree-preserving, and that $C_{\underline{n}}(D)$ is an object in the category Kom(Mat($\Kof_{/h}$)) whose objects are formal direct sums of objects in $\Kof_{/h}$. That is, $C_{\underline{n}}(D)$ is a complex of (direct sums of) formal complexes in $\Kof_{/h}$.

\begin{theorem}(cf.~\cite[Theorem 1]{BW})\label{thm:link inv}
The isomorphism class of the cochain complex $C_{\underline{n}}(D)$ is an invariant of the colored framed link $(L, \underline{n})$.
\end{theorem}

\noindent
\textit{Proof.} Let $D$ and $D'$ be diagrams representing isotopic colored framed links. 
Then, for any dot-row vector $\underline{s}$, the cable diagrams $D{\underline{s}}$ and $D'{\underline{s}}$ represent isotopic links, thus the formal complexes $[D{\underline{s}}]$ and $[D'{\underline{s}}]$ constructed using the $sl(2)$ cohomology theory are isomorphic as objects in $\Kof_{/h}$. 
The isotopy between the links represented by $D_{\underline{s}}$ and $D'_{\underline{s}}$ induces an isotopy between the cobordisms appearing in the definition of the differentials of $C_{\underline{n}}(D)$ and $C_{\underline{n}}(D')$.
Thus, complexes $C_{\underline{n}}(D)$ and $C_{\underline{n}}(D')$ are isomorphic. \hfill$\square$
\bigskip

To obtain a cohomology theory and a computable invariant we apply a functor to switch from the geometric category to an algebraic one. Specifically, we  apply the functor $\calF \co \foams_{/\ell} \to \rmod$ used in the $sl(2)$ foam cohomology theory. Denote by $\calF C_{\underline{n}}(D)$ the resulting complex, and by  $H_{\underline{n}} (D): = H(\calF C_{\underline{n}}(D))$ its cohomology. The cochain objects of the complex $\calF C_{\underline{n}}(D)$ are doubly-graded $R$-modules, and therefore, $H_{\underline{n}} (D)$ is a triply-graded $R$-module
\[ H_{\underline{n}} (D) = \bigoplus_{i,j,k \in \Z} H^{i,j,k}(D), \]
 where $i$ is the cohomological degree of $H_{\underline{n}} (D)$, and $j$ and $k$ are the cohomological and polynomial degrees, respectively, of the cochain objects of $\calF C_{\underline{n}}(D)$. Using Theorem~\ref{thm:link inv} and the fact that the functor $\calF$ is degree-preserving, we obtain that the isomorphism class of $H_{\underline{n}} (D)$ is independent on the diagram $D$ of the framed link $L$, and that is an invariant of $(L, \underline{n})$.

We form a three variable polynomial
\[P_{(L, \underline{n})}(r, t, q): = \sum_{i,j,k}r^it^jq^k \rk(H^{i,j,k}(D)),\]
and define the total graded Euler characteristic of $\calF C_{\underline{n}}(D)$ by
 \[ \chi(\calF C_{\underline{n}}(D)): = \sum_{i,j,k}(-1)^{i+j}q^k \rk(H^{i,j,k}(D)).\]
Then we have that $\chi(\calF C_{\underline{n}}(D)) = P_{(L, \underline{n})}(-1, -1, q)$. Moreover, $P_{(L, \underline{n})}(-1, -1, q) = J_{\underline{n}}(L)$, as shown below.

\begin{corollary}
The Euler characteristic of $\calF C_{\underline{n}}(D)$ is the colored Jones polynomial $J_{\underline{n}}(L)$.
\end{corollary}

\noindent
\textit{Proof.}    
\begin{eqnarray*}
\chi(\calF C_{\underline{n}}(D)) &=& \sum_{i,j,k}(-1)^{i+j}q^k \rk(H^{i,j,k}(D))\\
&=& \sum_{i}(-1)^i\sum_{\underline{s},\,\, |\underline{p}(\underline{s})| = i} \chi(\calF [\, D_{\underline{s}}\, ])\\
&=& \sum_{i}(-1)^i\sum_{\underline{s},\,\, |\underline{p}(\underline{s})| = i} \chi(\calF [\, D^{\underline{n} - 2 \underline{p}(\underline{s})}\, ])\\
&=& \sum_{\underline{k}} (-1)^{|\underline{k}|}\binom {\underline{n} - \underline{k}} {\underline{k}} \chi(\calF [\, D^{\underline{n} - 2 \underline{k}}\, ])\\
&=& \sum_{\underline{k} = \underline{0}}^{\left \lfloor  \underline{n}/2 \right \rfloor} (-1)^{|\underline{k}|}\binom{\underline{n} - \underline{k}} {\underline{k}} J(L^{\underline{n}-2\underline{k}}),
\end{eqnarray*}
where $\underline{k} = (k_1, k_2, \dots, k_l)$ is a vector whose entries are non-negative integers, and $|\underline{k}| = \sum_i k_i$.
Therefore, we obtain that $\chi(\calF C_{\underline{n}}(D)) = J_{\underline{n}}(L)$.
\hfill$\square$


\section{Colored tangle cohomology} \label{sec:colored tangle-coho}

\subsection{The case of tangles without closed components} \label{ssec:circle-free tangles}

In this section we consider oriented framed tangles $T$ without closed components (unless $T$ is a knot itself) whose strands are colored by the same natural number $n$. The $sl(2)$ foam cohomology theory is a ``local'' theory in the sense that is defined for arbitrary tangles, therefore it can be used to imitate the construction in Section~\ref{sec:colored link-coho} and associate to a diagram $D$ of a colored oriented framed tangle $(T,n)$ a complex $C_n(D)$ of formal complexes, and then a complex $\calF C_n(D)$ of doubly-graded $R$-modules.

\begin{figure}[ht]
\begin{center}
\includegraphics[height=4.2in, angle=90]{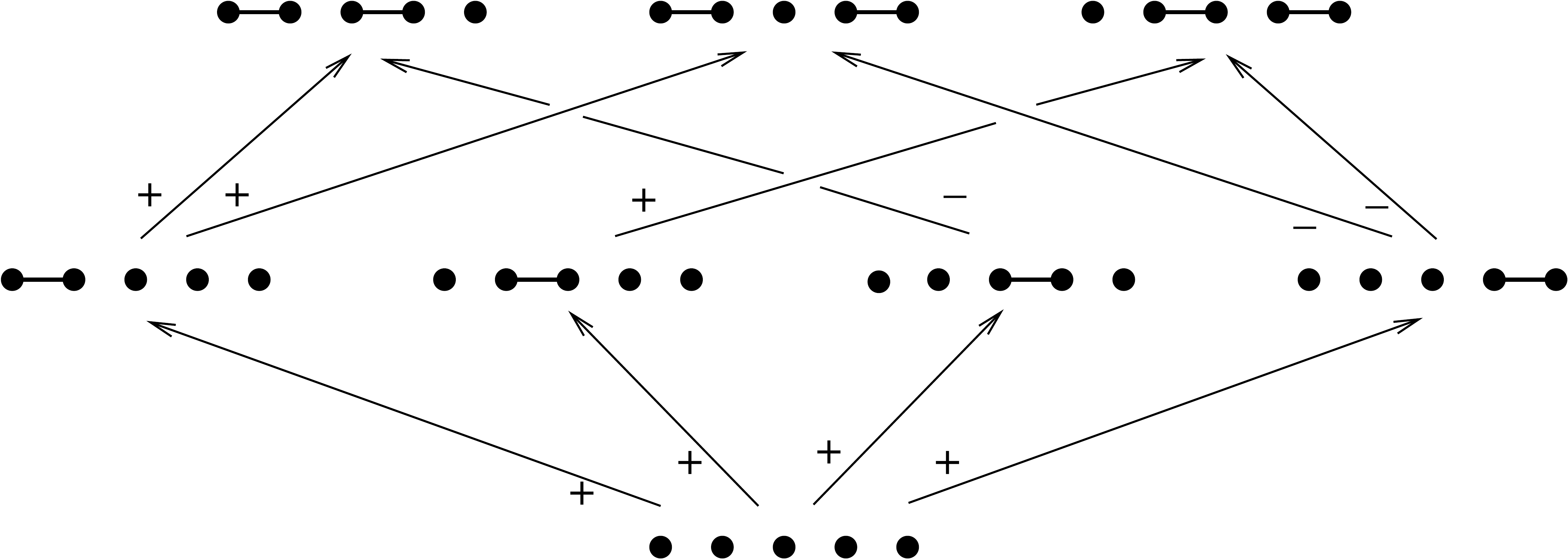} \\
\end{center}
\caption{The graph $\Gamma_5$}\label{fig:graph mono-colored tangle}
\end{figure}

Consider the graph $\Gamma_n$ whose vertices are marked by all dot-rows $s$ corresponding to $n$ (that is, dot-row vectors $\underline{s} = (s)$ with one component). Figure~\ref{fig:graph mono-colored tangle} displays the graph $\Gamma_5$. 

Associate the cable-diagram $D_s = D^{n-p(s)}$ to a dot-row $s$ in $\Gamma_n$. Diagram $D_s$ is the $(n-p(s))$-parallel cable of $D$, where there is a parallel cable-strand for each unpaired dot in $s$. Strand 1 is oriented in the same way as $D$, strand 2 is oppositely oriented, strand 3 is oriented as strand 1, etc. Below we show a tangle diagram $D$ and its 3-parallel cable $D^3$.

\[D = \qquad \raisebox{-25pt} {\includegraphics[height=0.8in]{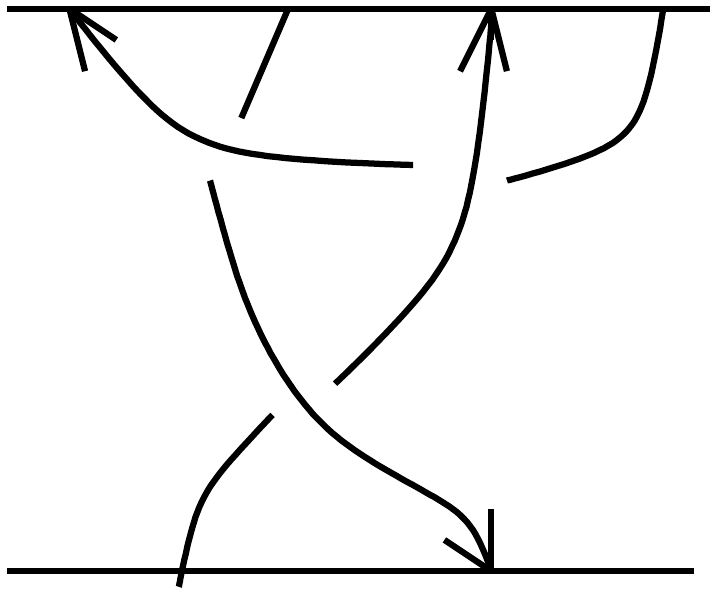}} \hspace{1in} D^3 = \raisebox{-25pt}{\includegraphics[height=0.8in]{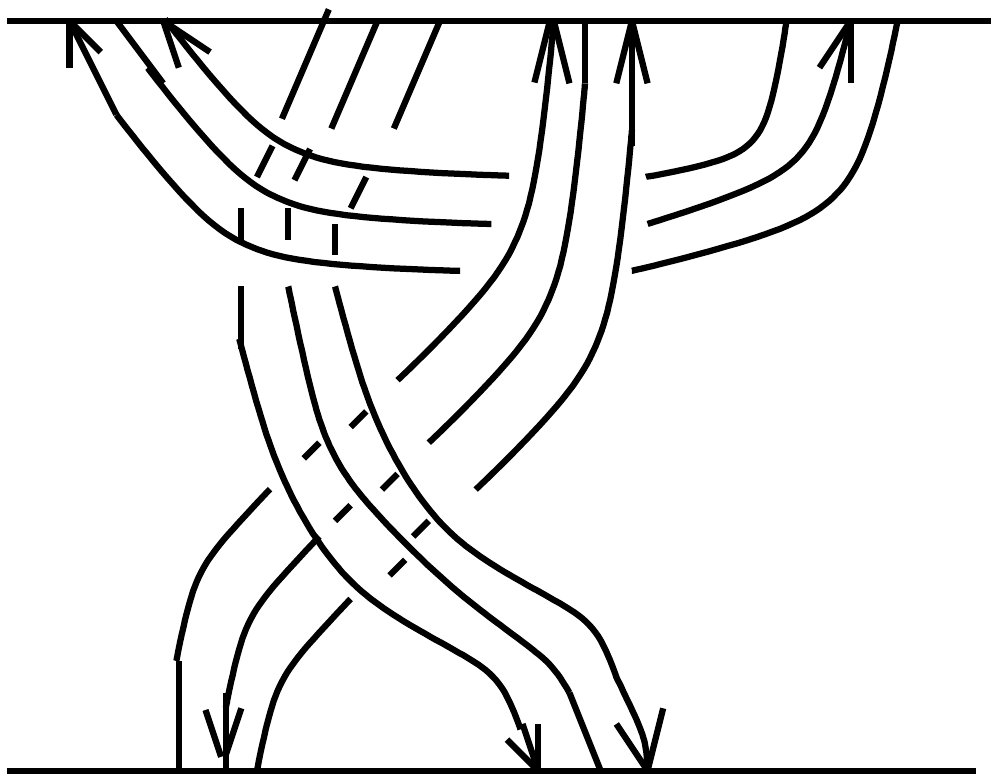}} \qquad\]

\begin{figure}[ht]
\begin{center}
\includegraphics[height=1.5in]{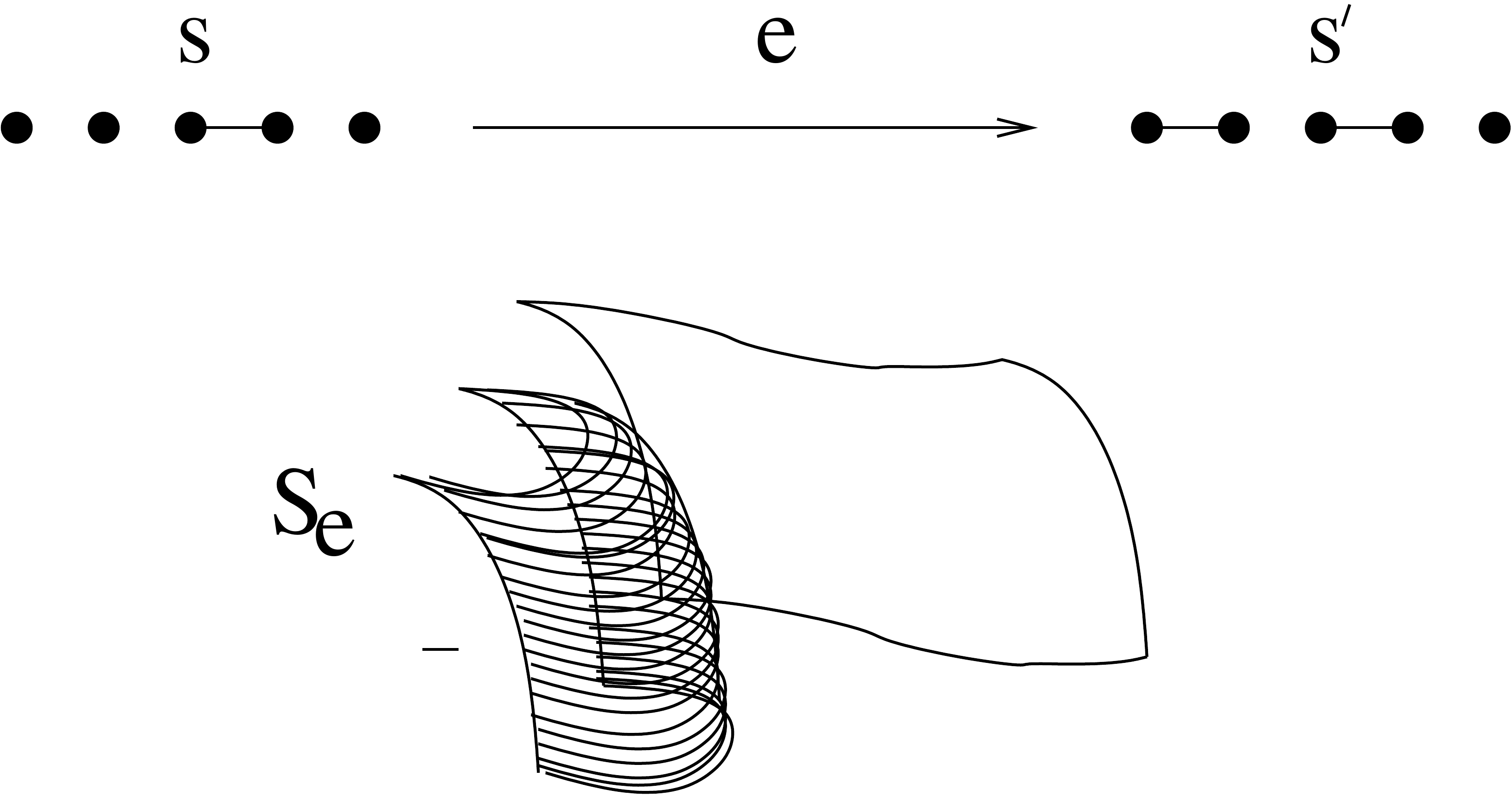}
 \end{center}
 \caption{The map $S_e$}\label{fig:edge-map}
 \end{figure}

For an edge $e \co s \to s'$, denote by $o(s,s')$ the number of pairs in $s$ to the right of the only pair in $s' \setminus s$. The map $S_e$ associated to an oriented edge $e \co s \to s'$ in the graph $\Gamma_n$ is a tangle cobordism from $D_s$ to $D_{s'}$, multiplied by $o(s,s')$. This cobordism is the identity everywhere except at the two adjacent strands in $D_s$ corresponding to the only pair in $s' \setminus s$, where the map is the cobordism obtained by contracting the two strands. In Figure~\ref{fig:edge-map} we show such a map $S_e$ for the rather boring $(1,1)$-identity tangle colored by $n=5$.

 To the latter graph we apply now the morphism $[ \, \cdot  \,]$ and form the complex $C_{n}(D)$ for the colored tangle theory. The cochain objects are given by
  \[ C^i_{n}(D): =  \bigoplus_{s,\, p(s) = i} [D_{s}].\]
The map $d^i : C^i_{n}(D) \to C^{i+1}_{n}(D)$ is a formal sum of all morphisms $[S_e]$ corresponding to edges $e$ at height $i$, where $[S_e] \co [D_s] \to [D_{s'}]$.

 \begin{proposition}\label{thm:tangleinv}
 The complex $C_n(D)$ is an invariant of the colored framed tangle $(T, n)$, up to homotopy. That is, if $D$ and $D'$ are isotopic colored framed tangle diagrams, then the cochain complexes $C_n(D)$ and $C_n(D')$ are isomorphic as objects in the category $\mbox{Kom}_{/h}(\mbox{Mat}(\Kof_{/h}))$.
 \end{proposition}
 
 \noindent
 \textit{Proof.} The proof is exactly the same as that of Theorem~\ref{thm:link inv}, thus we omit it. \hfill $\square$
 \bigskip

Finally, the colored tangle cohomology is obtained by applying the functor $\calF$ to the geometric invariant $C_n(D)$ of $T$. This yields a complex $\calF C_n(D)$ whose cochain objects are doubly-graded $R$-modules, and we take its cohomology. 
\begin{corollary}
The isomorphism class of the cohomology group $H_n(D): = H(\calF C_n(D))$ is a triply-graded invariant of the colored framed tangle $(T, n)$.
\end{corollary}

\emph{Remark}: If the tangle $T$ is a knot $K$ then the invariants $C_n(D)$ and $\calF C_n(D)$ agree with their analogues constructed in Section~\ref{sec:colored link-coho} for the colored link $(K, n)$ with one component.

\begin{corollary}
If the tangle $T$ is a knot $K$, then the graded Euler characteristic of $\calF C_n(D)$ is the colored Jones polynomial $J_n(K)$.
\end{corollary}

\subsection{Behavior under tangle composition}\label{ssec:behavior}

The goal of this section is to show that the geometric colored invariant of tangles defined in Section~\ref{ssec:circle-free tangles} composes well under (vertical) tangle composition. Specifically, let $D_1$ and $D_2$ be tangle diagrams corresponding to colored oriented framed tangles $(T_1, n)$ and $(T_2, n)$. Moreover, suppose that the vertical composition $D_1 \circ D_2$ is defined: 
\[ \raisebox{-6pt}{\includegraphics[height=0.25in]{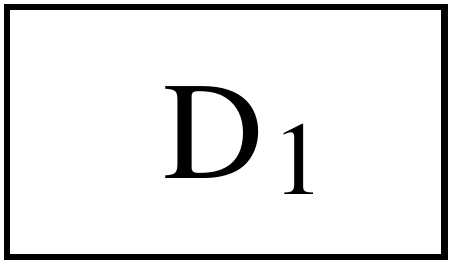}} \circ \raisebox{-6pt}{\includegraphics[height=0.25in]{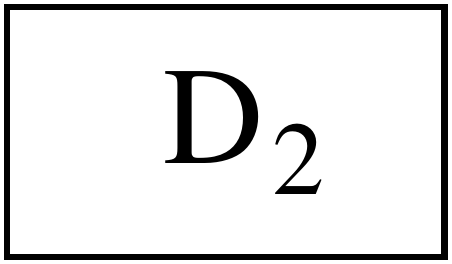}} = \raisebox{-15pt}{\includegraphics[height=0.45in]{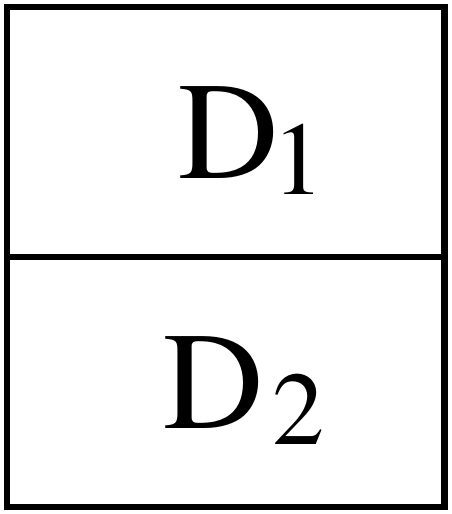}} \]

Here we consider again tangles with no closed components, therefore we need to impose that not only $D_1$ and $D_2$ are free of closed components, but also the resulting diagram $D_1 \circ D_2$. 

We show that there exists a binary operation $*$ defined on $\mbox{Kom}_{/h}(\mbox{Mat}(\Kof_{/h}))$ such that $C_n(D_1) * C_n(D_2) = C_n(D_1 \circ D_2)$, whenever the vertical composition $D_1 \circ D_2$ is defined . This operation is defined as follows.
 Let $C_n(D_1) = (C_n^i(D_1), d_1^i)$ and $C_n(D_2) = (C_n^i(D_2), d_2^i)$ where
\[C_n^i(D_1) = \displaystyle \bigoplus_{s, p(s) = i} [D_{1,s}] \,\,\, \mbox{and} \,\,\, C_n^i(D_2) = \displaystyle \bigoplus_{s, p(s) = i} [D_{2,s}]\] and \[ d_1^i (v_1) = \sum_{e} [S_{1,e}](v_1), \,\,\, d_2^i (v_2) = \sum_{e} [S_{2,e}](v_2)\,\,\,\, \mbox{for}\,\, v_1\in [D_{1,s}], \,\, v_2\in [D_{2,s}]\]
where the formal sum above is taken over all edges $e$ with tail $s$. 

Denote by
\[\calC^i = \displaystyle \bigoplus_{s, p(s) = i}([D_{1,s}] \otimes_R [D_{2,s}])\]
and
 \[\phi^i (v_1 \otimes v_2): =  \sum_{e} [S_{1,e}](v_1) \otimes_R [S_{2,e}] (v_2), \,\,\,\ \mbox{for}\,\, v_1 \otimes v_2 \in [D_{1,s}]  \otimes_R [D_{2,s}]  \]
 where the sum is taken over all edges $e$ with tail $s$, and define $C_n(D_1) * C_n(D_2): = (\calC^i , \phi^i)$.

Here $v_1$ and $v_2$ are resolutions---web diagrams---of $D_{1,s}$ and $D_{2,s}$, respectively, and $v_1 \otimes v_2$ stands for the resolution of the $D_{1,s} \circ D_{2,s}$ obtained by gluing (composing vertically) the webs $v_1$ and $v_2$ along their common boundary. The $k$-th direct summand of $\calC^i$, $[D_{1,s}] \otimes_R [D_{2,s}]$, is the \textit{ formal tensor product} of the $k$-th direct summands of $C_n^i(D_1)$ and $C_n^i(D_2)$. Specifically, the operation $\otimes$ here is the ``tensor product" operation induced on formal complexes by the composition operation on the canopoly $\foams_{/\ell}$, and which follows the gluing pattern used to make the tangle diagram $D_1 \circ D_2$ from $D_1$ and $D_2$ (this tensor product operation is possible by Proposition~\ref{prop:planar-algebra}.)

Let $e \co s \to s'$ be an edge, with associated cobordisms $S_{1,e} \co D_{1,s} \to D_{1,s'}$ and $S_{2,e} \co D_{2,s} \to D_{2,s'}$, and their induced maps $[S_{1,e}] \co [D_{1,s}] \to [D_{1,s'}]$ and $[S_{2,e}] \co [D_{2,s}] \to [D_{2,s'}]$, respectively. By Proposition~\ref{prop:canopoly}, we have that the tensor product of these maps 
\[[S_{1,e}] \otimes [S_{2,e}] \co [D_{1,s}] \otimes_R [D_{2,s}] \to  [D_{1,s'}] \otimes_R [D_{2,s'}]\] is equal to the following map, up to homotopy
\[ [S_{1,e} \circ S_{2,e}] \co [D_{1,s} \circ D_{2,s}] \to [D_{1,s'} \circ D_{2,s'}],\] 
where $S_{1,e} \circ S_{2,e}$ is the cobordism obtained by ``gluing'' $S_{1,e}$ and $S_{2,e}$ following the same pattern used to ``glue" (compose) the tangle diagrams which are the source and target of these cobordisms. Therefore, 
\[[S_{1,e}](v_1) \otimes_R [S_{2,e}] (v_2) = [S_{1,e} \circ S_{2,e}] (v_1 \otimes v_2).\]

We remark that $*$ is the ``direct sum" operation induced on complexes in Kom(Mat($\Kof_{/h}$)) by the composition operations on canopolies $\Kof_{/h}$ and $\cob^4_{/i}$. The following result holds at once.

\begin{proposition}
$C_n(D_1) * C_n(D_2)$ is a cochain complex. Moreover, $C_n(D_1) * C_n(D_2) = C_n(D_1 \circ D_2)$ up to homotopy.
\end{proposition}

\emph{Remark}:
We showed that the geometric colored invariant $C_n(T)$ of mono-colored framed tangles $T$ with no closed components has good composition properties, therefore it is suitable for efficient calculations of the colored cohomology groups of a framed knot. Specifically, we cut a colored oriented framed knot $(K, n)$ (using horizontal lines) into subtangles $(T_i, n)$, compute the invariants $C_n(T_i)$ and assemble them into $C_n(K)$, as prescribed in this subsection. Before the assembling operation, we simplify each $C_n(T_i)$ as much as possible by simplifying the cochain objects of  $C_n(T_i)$ (thus we simplify the formal complexes $[D_{i,s}] \in \Kof$ by making use of the ``delooping'' and ``Gaussian elimination'' procedures, as described in~\cite{C3}). Once that is taken care of, we apply the functor $\calF$ to arrive at the complex $\calF C_n(D)$, and take its cohomology.

\subsection{The case of arbitrary tangles} \label{ssec:general case}

In this section we show that we can do for links what we did for knots, namely we show that the ``divide and conquer" approach can be used to compute more efficiently the colored cohomology groups of an oriented framed link $(L, \underline{n})$. For that, we need a colored theory for arbitrary tangles, that is, tangles that might have closed components and whose strands might be colored by distinct natural numbers.

Let $L = T_1 \circ \ldots \circ T_k$ be an oriented framed link with $l$ components, regarded as a vertical composition of $k$ tangles.  Number the components of $L$ from $1$ to $l$, and color the $i$-th component by $n_i \in \N$.  Denote the colored link by $(L, \underline{n})$, where $\underline{n} = (n_1,  \ldots, n_l)$. 

The arcs of a subtangle $T_j$ correspond to certain link components, and thus receive the induced color from $L$. Denote by $\underline{n}_j$ the coloring of $T_j$ induced from the coloring of $L$, and denote the resulting colored tangle by $(T_j, \underline{n}_j)$.  Whenever a link component $K_i$ has representative arcs in $T_j$, we say that $K_i$ is \textit{represented} in $T_j$. The vector $\underline{n}_j$ has $l$ entries, and if all components of $L$ are represented in $T_j$, then $\underline{n}_j = \underline{n}$. Otherwise, if some link component $K_i$ is not represented in $T_j$ then the $i$-th entry of $\underline{n}_j$ is $0$, and all the other entries agree with the corresponding entries in $\underline{n}$. Thus $\underline{n}_j \in (\N \cup \{0\})^l$ for all $1 \leq j \leq k$.

\begin{figure}[ht]
\begin{center}
\includegraphics[height=1.8in]{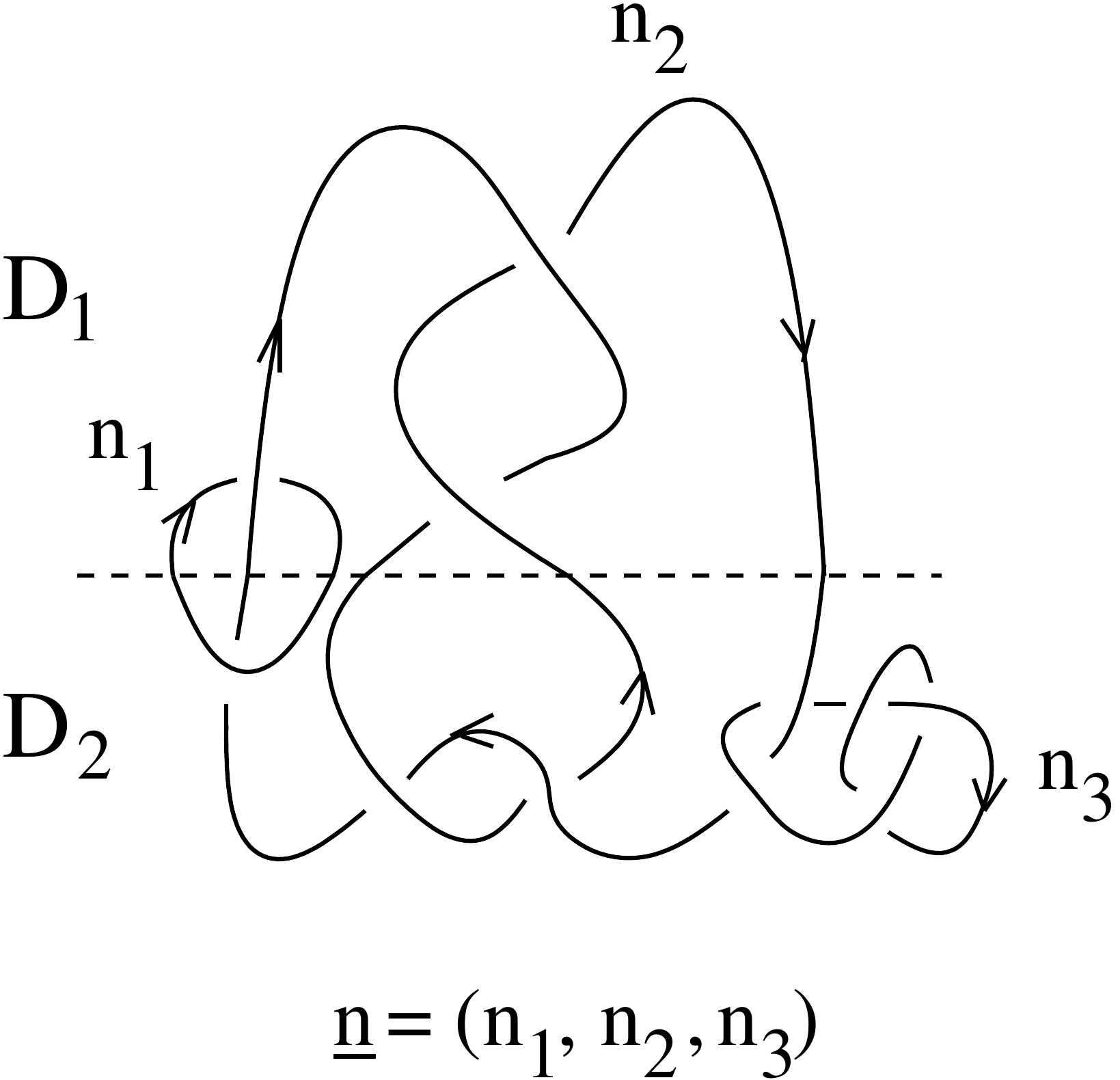} \quad \includegraphics[height=1.4in]{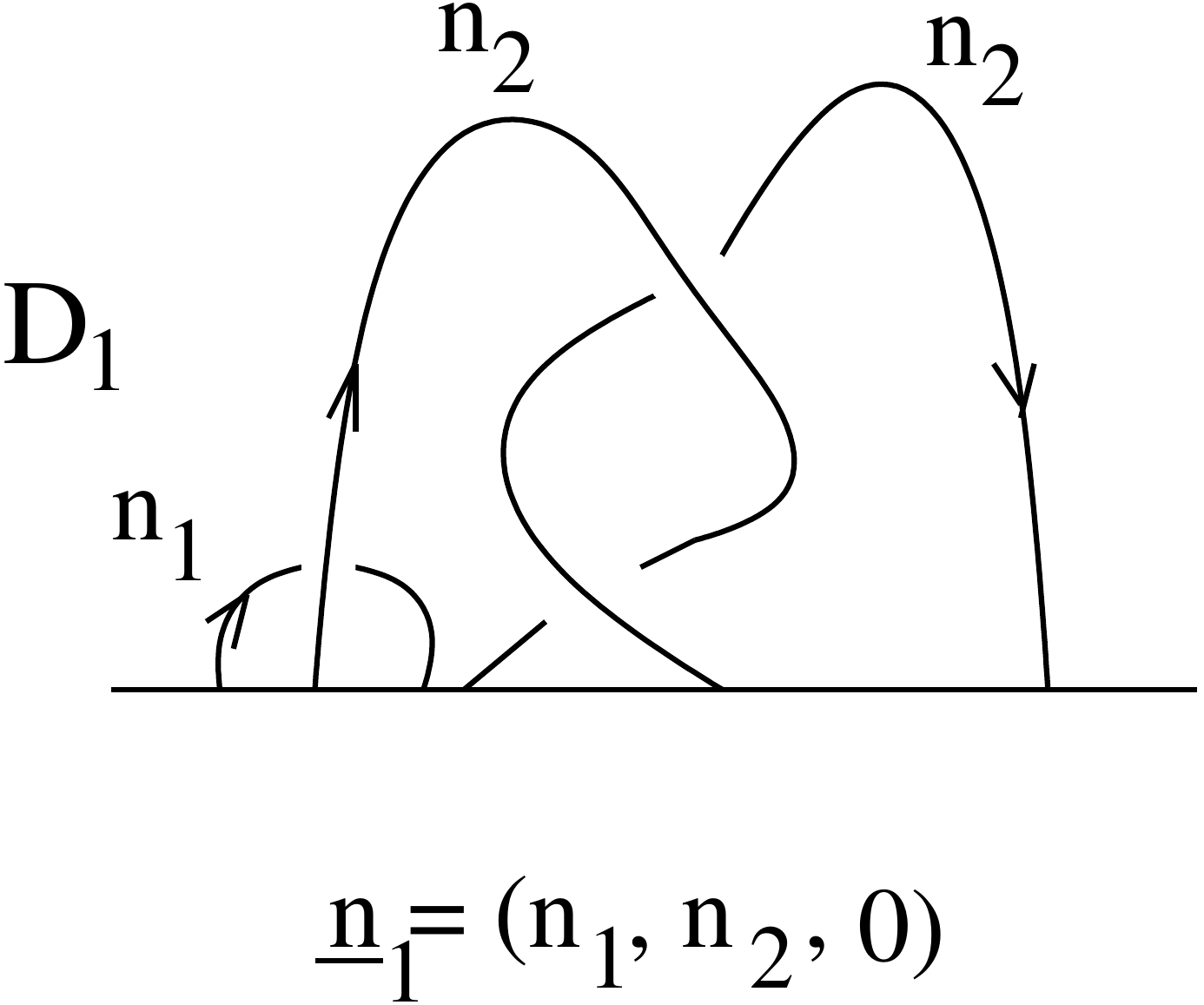} \quad  \includegraphics[height=1.1in]{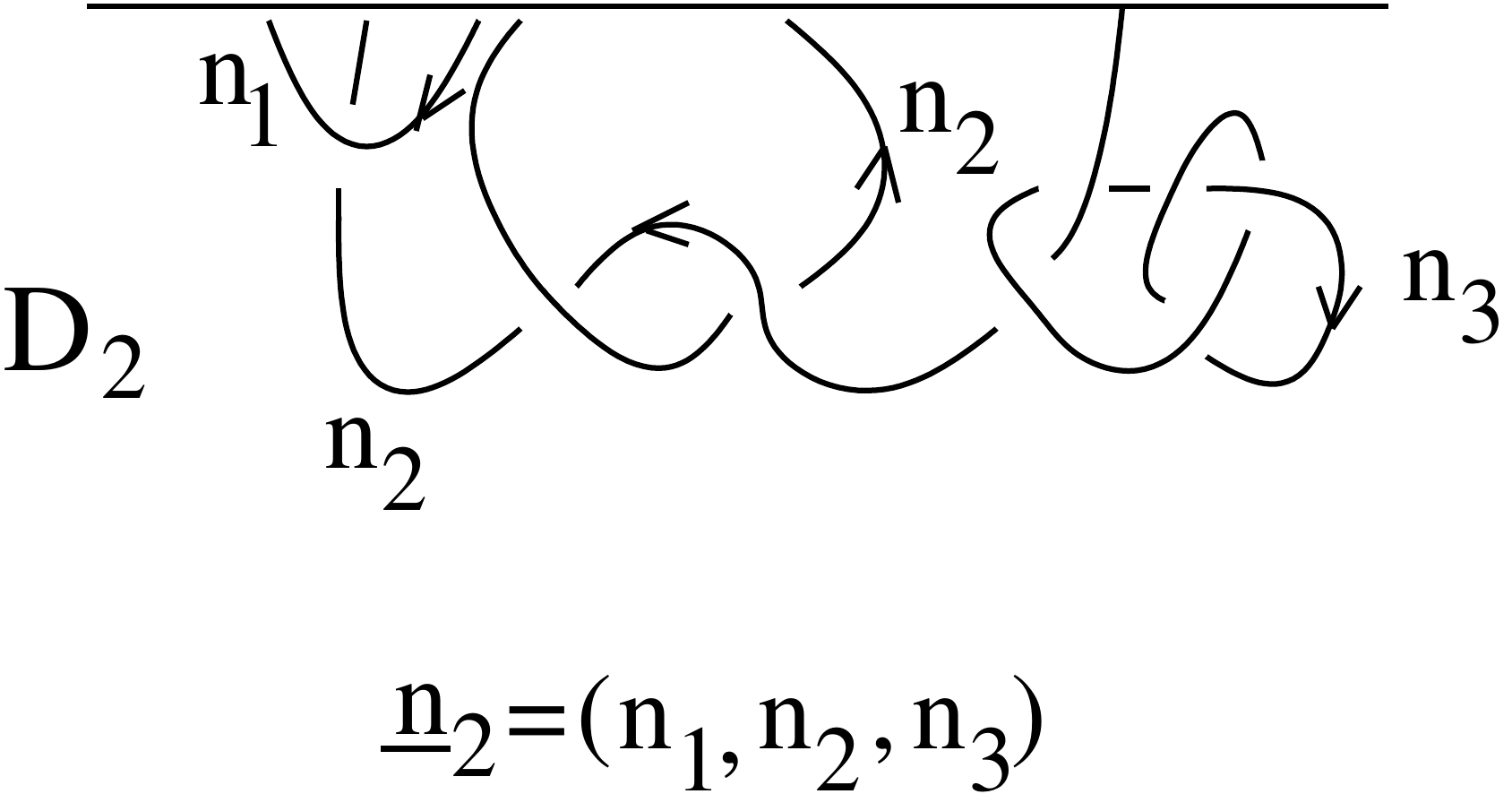} 
\end{center}
\caption{A link decomposed into subtangles}\label{fig:link}
\end{figure}

Let $D$ be a diagram of $L$ whose blackboard framing corresponds to the given framing of $L$, and write $D = D_1 \circ  \dots \circ D_k$, where $D_j$ is a diagram of $T_j$. In Figure~\ref{fig:link} we show a colored link diagram with three components, decomposed into two colored tangle diagrams.
 
 As before, let $\underline{s} = (s_1,  \dots, s_l)$ be a dot-row vector containing dot-rows $s_i$ with $n_i$ dots. Let $\underline{p}(\underline{s}) = (p(s_1),  \dots, p(s_l))$, where $p(s_i)$ is the number of pairs in $s_i$, and denote by $|\underline{p}(\underline{s})| = p(s_1) +  \ldots + p(s_l)$.
 
 For each diagram $D_j$ consider the (same) graph $\Gamma_{\underline{n}}$ whose vertices are marked by all dot-row vectors $\underline{s}$ corresponding to $\underline{n}$ (as constructed in Section~\ref{sec:colored link-coho}). Having a common graph for all subtangle diagrams $D_j$ is essential for obtaining a well-defined composition operation of the geometric colored invariants of (arbitrary) tangles, and therefore, in recovering the colored Jones polynomial of a link. 

Consider the tangle diagram $D_j$, where $j$ is fixed. To a dot-row vector $\underline{s} \in \Gamma_n$ associate the cable diagram $D_{j, \underline{s}}: = D_j^{\underline{n}_j - 2\underline{p}(\underline{s})}$ formed by taking the $(n_i - 2p(s_i))$-parallel cable of each arc in  $D_j$ colored by $n_i$. If the $i$-th entry in $\underline{n}_j$ is $0$ (that is, if the component $K_i$ of $L$ is not represented in $D_j$), or equivalently, if there are no arcs in $D_j$ colored by $n_i$, then take the $(n_i - 2p(s_i))$-cable of the \textit{empty} tangle diagram for the missing arcs.
To an edge $e \co \underline{s} \to \underline{s}'$ associate the cobordism $S_{j,e}: D_{j, \underline{s}} \to D_{j, \underline{s}'}$ given by contracting the neighboring strands in $D_{j, \underline{s}}$ corresponding to the pair in $\underline{s}' $ but not in $\underline{s}$.
Finally, multiply each cobordism $S_{j,e}$ by $(-1)^{o(\underline{s}, \underline{s}')}$, where $o(\underline{s},\underline{s}')$ is the number of pairs in $\underline{s}$ to the right and above of the only pair in $\underline{s}' \setminus \underline{s}$ (see Figure~\ref{fig:graph for a link}). 

We are ready now to form the complex $C_{\underline{n}_j}(D_j)$ by applying the morphism $[ \, \cdot  \,]$ of $sl(2)$ foam cohomology theory. Let $C_{\underline{n}_j}(D_j)  = (C^i_{\underline{n}_j}(D_j), d^i_j )$. Then
\[ C^i_{\underline{n}_j}(D_j): =  \bigoplus_{\underline{s}} [D_{j, \underline{s}}]\]
where the sum is taken over all dot-row vectors $\underline{s}$ with $|\underline{p}(\underline{s})| = i$, and the map $d^i_j : C^i_{\underline{n}_j}(D_j) \to C^{i+1}_{\underline{n}_j}(D_j)$ is a formal sum of all morphisms $[S_{j,e}]$ corresponding to edges $e$ with tail $\underline{s}$. 

The following proposition is proved much as Theorem~\ref{thm:link inv}.

\begin{proposition}\label{thm:arbitrary tangle inv}
$C_{\underline{n}_j}(D_j)$ is a complex whose isomorphism class is an invariant of the colored framed tangle $(T_j, \underline{n}_j)$.
\end{proposition}

The operation $*$ of $C_{\underline{n}_j}(D_j)$ and $C_{\underline{n}_{j+1}}(D_{j+1})$, for $1 \leq j \leq k-1$ goes as follows. Denote by  $(\calC^i, \phi^i) = C_{\underline{n}_j}(D_j) * C_{\underline{n}_{j+1}}(D_{j+1})$, and define
\[\calC^i : = \displaystyle \bigoplus_{\underline{s}}([D_{j, \underline{s}}] \otimes_R [D_{j+1, \underline{s}}])\]
where the sum is taken over all dot-row vectors $\underline{s}$ such that $| \underline{p}(\underline{s})| = i$. Moreover, consider the map $\phi^i \co \calC^i \to \calC^{i + 1}$ which is the formal sum of all morphisms $[S_{j,e}] \otimes [S_{j + 1,e}]$ corresponding to all edges $e$ with tail $\underline{s}$. It follows that $(\calC^i, \phi^i)$ is a complex whose construction is modeled by the \textit{formal direct sum} of complexes $C_{\underline{n}_j}(D_j)$ and $C_{\underline{n}_{j+1}}(D_{j+1})$, as explained in Section~\ref{ssec:behavior}.

Moreover, the isomorphism class of $C_{\underline{n}_j}(D_j) * C_{\underline{n}_{j+1}}(D_{j+1})$ is the colored invariant of $T_j \circ T_{j + 1}$, by construction. Putting all together, $C_{\underline{n}_1}(D_1) * \ldots * C_{\underline{n}_{k}}(D_{k}) = C_{\underline{n}}(D)$ up to homotopy, and therefore, $C_{\underline{n}_1}(D_1) * \ldots * C_{\underline{n}_{k}}(D_{k})$ is an up-to-homotopy invariant of the colored link $(L, \underline{n})$.

It is important to remark that, as in the case of knots, we simplify as much as possible each of the invariants $C_{\underline{n}_j}(D_j)$, $1\leq j \leq k$, before we assemble  them into the invariant of the colored link $(L, \underline{n})$. Applying the functor $\calF$ we obtain a complex $\calF (C_{\underline{n}_1}(D_1) * \ldots * C_{\underline{n}_{k}}(D_{k}))$ of doubly-graded $R$-modules and homomorphism between them, and we can take its cohomology. The isomorphism class of $H(\calF (C_{\underline{n}_1}(D_1) * \ldots * C_{\underline{n}_{k}}(D_{k})))$ is an invariant of the framed colored link $(L, \underline{n})$, and its total graded Euler characteristic is the colored Jones polynomial $J_{\underline{n}}(L)$.
\bigskip

\textbf{Acknowledgements.} The author gratefully acknowledges NSF support  through grant DMS 0906401. She would also like to thank A. Shumakovitch for suggesting to consider working with tangles (as opposed  to just knots and links) for the sake of obtaining efficient computations, and to the referee for providing valuable comments.


\begin{flushleft}
Department of Mathematics, 
California State University, Fresno, CA 93740, USA \\
E-mail address: {\tt ccaprau@csufresno.edu} 
\end{flushleft}

\end{document}